\documentclass[11pt]{article}

\usepackage{mathrsfs}
\usepackage{latexsym}
\usepackage{amsmath,amsthm}
\usepackage{amsfonts}
\usepackage{url}

\usepackage{amssymb}
\usepackage{mytheo}
\newcommand{\etab}{\bar\eta}

\usepackage{bbm}
\newcommand{\chr}{\boldsymbol{\mathbbm{1}}} % characteristic function
\newcommand{\pred}[1]{\chr_{\left\{ #1 \right\}}}
\newcommand{\prs}{\vec{P}}
\newcommand{\pr}[1]{\prs\!\tlprn{#1}}

\renewcommand{\P}{\prs}

\newcommand{\trn}{^{\!\mathsf{T}}} %operator transpose
\newcommand{\inv}{^{-1}} %inverse

\newcommand{\nn}[1]{[#1]}
\newcommand{\ddel}[2]{\frac{\del#1}{\del#2}}
\newcommand{\grad}{\nabla}
\newcommand{\del}{\partial}

\newcommand{\TV}[1]{\nrm{#1}_{\textrm{{\tiny \textup{TV}}}}}

\newcommand{\Lip}[1]{\nrm{#1}_{\textrm{{\tiny \textup{Lip}}}}}

\renewcommand{\omit}[2]{ \bar{#1}_{#2} }

\usepackage{geometry}
\makeatletter

\newcommand{\ben}{\begin{enumerate}}
\newcommand{\een}{\end{enumerate}}
\newcommand{\bit}{\begin{itemize}}
\newcommand{\eit}{\end{itemize}}

\newcommand{\labs}{\left| \vphantom{\sum_a^b} \right.}
\newcommand{\rabs}{\left. \vphantom{\sum_a^b} \right|}
\newcommand{\lparen}{\left( \vphantom{\sum} \right.}
\newcommand{\rparen}{\left. \vphantom{\sum} \right)}

\renewcommand{\vec}[1]{\bs{\mathrm{#1}}}

\newcommand{\basicspace}{\Omega}
\newcommand{\X}{\basicspace}
\newcommand{\supr}[1]{^{(#1)}}
\newcommand{\seq}[3]{(#1_{#2},\ldots,#1_{#3})}
\newcommand{\sseq}[3]{#1_{#2}^{#3}}  % short seq
\newcommand{\sumseq}[3]{\sum_{\sseq{#1}{#2}{#3}}}

\newcommand{\dsabs}[1]{\bigl| #1 \bigr|}
\newcommand{\nrm}[1]{\left\Vert #1 \right\Vert}
\newcommand{\dsnrm}[1]{\bigl\Vert #1 \bigr\Vert}

\newcommand{\iprod}[1]{\left\langle #1 \right\rangle}

\renewcommand{\b}[1]{\hat{#1}}

\newcommand{\calL}{\mathcal{L}}

\newcommand{\calX}{\mathcal{X}}
\newcommand{\calY}{\mathcal{Y}}

\newcommand{\tha}{\theta}
\newcommand{\R}{\mathbb{R}}
\newcommand{\N}{\mathbb{N}}

\newcommand{\beq}{\begin{eqnarray*}}
\newcommand{\eeq}{\end{eqnarray*}}
\newcommand{\beqn}{\begin{eqnarray}}
\newcommand{\eeqn}{\end{eqnarray}}
\newcommand{\paren}[1]{\left( #1 \right)}
\newcommand{\sqprn}[1]{\left[ #1 \right]}
\newcommand{\tlprn}[1]{\left\{ #1 \right\}}
\newcommand{\set}[1]{\tlprn{#1}}
\newcommand{\abs}[1]{\left| #1 \right|}

\newcommand{\floor}[1]{\ensuremath{\left\lfloor#1\right\rfloor}}
\newcommand{\gn}{\, | \,}

\newcommand{\ds}{\displaystyle}
\newcommand{\ts}{\textstyle}
\newcommand{\bs}{\boldsymbol}
\renewcommand{\th}{\ensuremath{^{\mathrm{th}}}~}

\newcommand{\hide}[1]{}
\newcommand{\oo}[1]{\frac{1}{#1}}

\newcommand{\one}{\text{i}}
\newcommand{\two}{\text{ii}}

\newcommand{\iia}{\ensuremath{\mathrm{(a)}}}
\newcommand{\iib}{\ensuremath{\mathrm{(b)}}}
\newcommand{\iic}{\ensuremath{\mathrm{(c)}}}
\newcommand{\iid}{\ensuremath{\mathrm{(d)}}}
\newcommand{\iie}{\ensuremath{\mathrm{(e)}}}
\newcommand{\iif}{\ensuremath{\mathrm{(f)}}}

\def\eps{\varepsilon}
\newcommand{\defeq}{=}
\newcommand{\tp}{\otimes}
\newcommand{\TP}{\bigotimes}

\title{Obtaining Measure Concentration from Markov Contraction}
\author{Aryeh Kontorovich\\
Department of Computer Science \\
Ben-Gurion University\\
Beer Sheva, Israel}

\begin{document}
\maketitle
\begin{abstract}
Concentration bounds for non-product, non-Haar measures are fairly recent: 
the first such result was obtained for contracting
Markov chains by Marton in 1996 via the coupling method. 
The work that followed, with few exceptions, also used coupling.
%Since then, several other such results have been proved; with few exceptions, these rely
%on coupling techniques.
Although this technique is of unquestionable utility as a theoretical tool, it 
is not always simple to apply.
%appears to have some limitations.
%Coupling has yet to be used to obtain bounds for more general Markov-type processes: hidden (or partially observed) Markov chains,
%Markov trees, etc. 
As an alternative to coupling, we use the elementary Markov contraction lemma 
%in a novel way, 
to obtain 
simple, useful, and apparently novel concentration results for various Markov-type processes. 
Our technique 
consists of expressing probabilities as matrix products and applying Markov contraction to these expressions;
thus it is fairly general and
%is generic and 
holds the potential to yield further results in this vein. 

\hide{
The problem of establishing measure concentration for dependent
processes under strong mixing conditions has been brought to
prominence by the pioneering work of Marton [[cites]] and others
[[cites]].  The two main approaches have been the coupling method and
bounds on contraction coefficients.  While the latter has been less
prominent recent work by [[cites]] have shown that it can deliver
powerful results, enabling us to control mixing coefficients which
are global properties of the process in terms of easily-calculated
local quantities.  In this paper, I unify these results on
establishing measure concentration by means of Markov contraction
coefficients, and extend them to various not-quite-Markov processes."
would be a big help to the reader

-technique basically consists of expressing probabilities as matrix products and applying Markov contraction to these expressions

-on some level, the distinction between our technique and the coupling method is semantic. By ``coupling'' we mean a 
construction of a joint distribution
-some ppl might refer to contraction as coupling
}
\end{abstract}

\section{Introduction}
\subsection{Background}
In 1996 Marton \cite{marton96}
published a concentration inequality for contracting Markov chains
 --- apparently, the first such result for a 
non-product, non-Haar measure. In the decade that followed, Marton and others
%a number of authors 
continued to 
%distill and expand 
deepen and broaden
a key
insight: analogues of the Azuma-Hoeffding-McDiarmid 
inequality 
\cite{azuma,hoeffding,mcdiarmid89}
for independent random variables
may be obtained for dependent ones, provided
a strong mixing condition holds.

To recall, the aforementioned inequality 
%states 
implies
that if $\mu$ is a product 
%measure
distribution
on $\X^n$ and
$f:\X^n\to\R$ satisfies
$\Lip{f}\leq 1$
under the 
%normalized 
Hamming metric, we have
\beqn
\label{eq:mcd}
\mu\set{\abs{f-\mu f}>t} &\leq& 2\exp(-2t^2/n).
\eeqn

In \cite{marton96}, Marton pioneered the transportation method for proving concentration inequalities. This technique
is in principle applicable to arbitrary nonproduct measures, and when applied to Markov chains 
$\mu$
with contraction coefficient $\tha<1$,
it yields
\beqn
\label{eq:marton}
\mu\set{\abs{f-M_f}>t} &\leq &
2\exp\sqprn{
-\frac{2}{n}\paren{
t(1-\tha)
-\sqrt{\frac{\log2}{2n}}}^2
},
\eeqn
where $M_f$ is a $\mu$-median of $f$. Since product 
distributions
%measures 
are
degenerate cases of Markov 
%measures 
chains
(with $\tha=0$), Marton's result
is a powerful generalization of (\ref{eq:mcd}).

The Markov contractivity condition $\tha<1$ implies strong mixing, and in a series of papers
\cite{marton98,marton03,marton04}, Marton gave other concentration results for dependent variables
under various metrics and types of mixing. In particular, Theorem 2 of \cite{marton98} gives a generic mixing condition
which implies a transportation inequality and therefore concentration.

Further progress in obtaining concentration from mixing was made, 
among others, in
\cite{chatterjee05,chazottes07,MR2511280,kontram06,kulske03,rio00,samson00}.
Using Stein's method for exchangeable pairs, 
Chatterjee 
\cite{chatterjee05}
obtained an elegant concentration inequality in terms of a Dobrushin-Shlosman type contractivity condition.
Samson \cite{samson00} was apparently the first to use explicit mixing coefficients in a concentration result.
Since these are central to this paper we define them without further delay; the (standard) notation is clarified in Section \ref{sec:notconv}.

Let $\mu$ be the joint distribution of 
$\seq{X}{1}{n}$, $X_i\in\X$.
For $1\leq i<j\leq n$
and $x\in\X^i$,
we denote by
$$\mu(\seq{X}{j}{n}\gn \seq{X}{1}{i}=x)
$$ 
the distribution
of $\seq{X}{j}{n}$ conditioned on $\seq{X}{1}{i}=x$. 
For $y\in\X^{i-1}$ and 
$w,w'\in\X$, 
define
\beq
%\nonumber
\eta_{ij}(y,w,w') &=&
\TV{
\mu(\seq{X}{j}{n}\gn \seq{X}{1}{i}={y w})-
\mu(\seq{X}{j}{n}\gn \seq{X}{1}{i}={y w'})
},
%\\&
\eeq
and
\beqn
\label{eq:etadef}
\etab_{ij} &=&
\sup_{y\in\X^{i-1},w,w'\in\X}
\eta_{ij}(y,w,w').
\eeqn
The coefficients $\etab_{ij}$, 
termed {\em $\eta$-mixing coefficients}\footnote{
That choice of terminology 
%may be 
is perhaps
suboptimal in light of the unrelated
notion of {\em $\eta$-weak dependence} of Doukhan et al. \cite{doukhan99},
but the sufficiently distinct contexts should prevent confusion.
}
in \cite{kontram06},
play a key role in several recent concentration results. 
Define $\Gamma$ and $\Delta$ to be upper-triangular $n\times n$ matrices,
with $\Gamma_{ii}=\Delta_{ii}=1$ and
\beq
\Gamma_{ij} = \sqrt{\etab_{ij}},
\qquad
\Delta_{ij} = \etab_{ij}
\eeq
for $1\leq i<j\leq n$.

In 2000, Samson \cite{samson00} proved that 
any distribution $\mu$ on $[0,1]^n$ and any convex
$f:[0,1]^n\to\R$ with $\Lip{f}\leq 1$ (with respect to $\ell_2$) satisfy
\beqn
\label{eq:samson}
\mu\set{\abs{f-\mu f}>t} &\leq& 2\exp\paren{-\frac{t^2}{2\nrm{\Gamma}_2^2}}
\eeqn
where $\nrm{\Gamma}_2$ is the $\ell_2$ operator norm.

In 2007, almost 
%contemporaneously 
synchronously
and using different techniques,
Chazottes et al. \cite{chazottes07} and
the author with K. Ramanan \cite{kontram06} showed that 
any distribution $\mu$ on $\X^n$ and any
$f:\X^n\to\R$ with $\sqrt n\Lip{f}\leq 1$ (with respect to the 
%normalized 
Hamming metric) satisfy
\beqn
\label{eq:kontram}
\mu\set{\abs{f-\mu f}>t} &\leq& 2\exp\paren{-\frac{t^2}{2\nrm{\Delta}_\infty^2}}
\eeqn
where $\nrm{\Delta}_\infty$ is the $\ell_\infty$ operator norm
($\nrm{\Delta}_\infty$ may be replaced by $\nrm{\Delta}_2$ and
%Chazottes et al. 
\cite{chazottes07} achieves
%obtain 
a better constant in the exponent).
More recently, \cite{MR2511280} obtained polynomial concentration bounds in terms of moments of the coupling time. The contraction condition may be cast more generally as a curvature 
\cite{MR2484937} or a ``metric ergodicity''
\cite{MR2683634}
condition. 
Some recent results on 
transportation methods in concentration 
include
\cite{MR2438906,djellout04,springerlink:10.1007/s00440-008-0159-5};
a survey
may be found in \cite{kolesnikov2010}.

The results (\ref{eq:samson}) and (\ref{eq:kontram}) are not readily comparable as they hold in different spaces for different 
metrics
with different normalization,
and the former requires convexity. They share the feature of establishing concentration for a wide class of measures, in terms
of the natural mixing coefficients $\etab_{ij}$. Indeed, since
\beqn
\label{eq:Delta-inf}
\nrm{\Delta}_\infty &=& \max_{1\leq i<n} (1 + \etab_{i,i} + \etab_{i,i+1} + \ldots + \etab_{i,n})
\eeqn
and
by the Ger\v{s}gorin disc theorem \cite{horn-johnson}
\beq
\nrm{\Gamma}_2^2 = \lambda_{\max}(\Gamma\trn\Gamma) &\leq&
\max_{1\leq i\leq n} \sum_{j=1}^n (\Gamma\trn\Gamma)_{ij},
\eeq
suitable upper estimates on $\etab_{ij}$ provide bounds for 
$\nrm{\Gamma}_2$
and
$\nrm{\Delta}_\infty$.

Aside from the straightforward observation (due to Samson) that the $\eta$-mixing coefficients are 
%controlled 
bounded
by the 
$\phi$-mixing ones
(see \cite{bradley05}), we are only aware of a few of cases where simple, readily computed estimates on $\etab_{ij}$ are given.
In particular, 
Samson \cite{samson00} controls $\etab_{ij}$ by the contraction coefficients of a Markov chain, and
Chazottes et al. \cite{chazottes07} give some estimates on $\etab_{ij}$ for various temperature regimes of Gibbs random fields.
The estimates quoted above are obtained via the coupling method --- which, while powerful, often requires some ingenuity to construct the 
requisite joint distribution, even in the simple case of a Markov chain \cite{marton96,samson00}. In some cases, the coupling may even
elude explicit construction \cite{chazottes07}.

As the random processes of interest become more complex, it becomes progressively more difficult to obtain estimates
on mixing coefficients via coupling. We are particularly interested in examining the $\eta$-mixing of several Markov-type
processes, motivated by statistical and computer science applications. Hidden Markov Models (HMMs) have been used in natural
language processing \cite{manning99,rabiner89} and signal processing \cite{Maybeck79}
for decades, with considerable success. 
Concentration bounds for Markov Chains (and more generally, HMMs) have
implications in machine learning and empirical process theory \cite{gamarnik03,kont07-thesis}.
A Markov-type process called the {\em Markov marginal process} (MMP)
in 
\cite{kont07-slln} 
underlies adaptive Markov Chain Monte Carlo simulations \cite{atchade05}; these evolve according to
an inhomogeneous Markov kernel, which in addition to time also depends on the path history. In a forthcoming work,
A. Brockwell and the author give strong laws of large numbers for MMPs in therms of the $\eta$-mixing coefficients.
Random processes indexed by trees have been attracting the attention of probability theorists for some time
\cite{benjamini94,pemantle95}, and the principal technical contribution of this paper is a bound on $\etab_{ij}$ for these
types of processes.

\hide{
Our results do not invoke the coupling method but rather rely on the Markov contraction lemma (Lemma \ref{lem:markcontr}). The
technique provides novel concentration bounds 
%for hidden and partially observed Markov chains as well as Markov trees -- 
for the processes listed above, and in some cases significantly sharpens
the estimates appearing in the literature.
% ---
%results which
%the coupling method has yet to yield or reproduce.
}
\begin{rem}
Although our results rely on
Markov contraction (Lemma \ref{lem:markcontr})
rather than on the coupling method,
on some level, the distinction 
%between our technique and the coupling method 
is semantic. 
%From conversations with experts it appears that
%what we call here ``Markov contraction''
%% (Lemma \ref{lem:markcontr}) 
%is commonly referred to as ``coupling''. 
The novelty of our method lies
in 
(i) avoiding any constructions (implicit or explicit) of joint distributions 
(ii) rewriting complicated sums as simple(r) matrix and tensor products
(iii) applying Lemma \ref{lem:markcontr} to the latter expressions.
Thus it seems that our method is sufficiently 
%distinct 
different
from classical coupling techniques
%, both in execution and results obtained, 
to merit
the terminological distinction.
\end{rem}
\hide{
-
By ``coupling'' we mean a 
construction of a joint distribution
-some ppl might refer to contraction as coupling
}

\subsection{Main results}
\label{subs:mainres}
%The main results of this paper are 
In this paper we present
estimates on the $\eta$-mixing coefficients $\etab_{ij}$ defined in (\ref{eq:etadef}), for
the various Markov-type processes mentioned above.
These bounds immediately imply concentration inequalities for a wide class of metrics and measures, via 
(\ref{eq:samson}) and (\ref{eq:kontram}).

The precise statements of the results 
%are technical 
require preliminary definitions
and are postponed until later sections.
The main technical contribution of this paper is Theorem \ref{thm:mtmain}, which bounds $\eta$-mixing coefficients for
Markov-tree processes, yielding what appears to be the first concentration of measure result for these.
However, 
we give equal priority to the 
goal of
%our desire is to highlight 
%highlighting 
presenting
Markov contraction as a 
%powerful 
versatile
new method
%a new method 
for bounding $\etab_{ij}$.
% as opposed to featuring any particular bound.
The nature of the bounds is to control
$\etab_{ij}$ --- a global function of the 
%measure 
distribution
$\mu$ --- by some local, easily computed contraction coefficients of $\mu$.
For example,
%Thus, 
let $\mu$ be an inhomogeneous Markov chain 
%given 
defined
by the transition kernels
$\set{p_i: 0\leq i< n}$,
which induces a 
%probability measure 
%density
distribution
on $\X^n$ by
\beq
\mu(x) &=& p_0(x_1)\prod_{i=1}^{n-1} p_i(x_{i+1}\gn x_i),
\qquad x\in\X^n.
\eeq
Define the  $i$\th 
{\em contraction coefficient}:
\beqn
\label{eq:thadef}
\tha_i &=&
\sup_{y,y'\in\X}\TV{p_i(\cdot\gn y)-p_i(\cdot\gn y')},
\qquad 1\leq i<n
.
\eeqn
This quantity turns out to control the $\eta$-mixing coefficients for
$\mu$:
\beq
\etab_{ij} &\leq& \tha_i\tha_{i+1}\ldots\tha_{j-1}
\eeq
--- a fact which is proved in \cite{samson00} using coupling. 
In \cite{kontram06} we gave
%We will give 
an (arguably simpler) alternative proof, which 
%will pave 
paves
the way for the several new results presented here.

%\subsection{Paper outline}
%\label{subs:outline}
This paper is organized as follows.
In Section \ref{sec:notconv} we summarise some basic notation used
throughout the paper. Some auxiliary lemmas are given in Section \ref{sec:contr}. The
remaining three sections deal with bounding $\etab_{ij}$ for
Markov chains (directed and undirected) and Markov tree processes.
%, and
%Markov marginal processes,
%%partially observed Markov chains
%respectively.

\subsection{Notation and definitions}
\label{sec:notconv}
Since the contribution of this paper is not measure-theoretic in nature, we henceforth take $\X$ to be a finite set.
Extensions to the countable case are quite straightforward \cite{kontram06} and the continuous case,
under mild assumptions, is not much more
difficult \cite{kont07-thesis,kont07-slln}.

We 
%will 
use the terms {\em measure} and {\em distribution} interchangeably; 
all measures are probabilities unless noted otherwise. If $\mu$ is
a measure on $\X^n$ and $f:\X^n\to\R$, we 
use the standard notation
\beq
\mu f &=& \int_{\X^n} f d\mu
\eeq
and
write
\beq
\mu\set{\abs{f-\mu f}>t}
\eeq
as a shorthand for
\beq
\mu\paren{\set{x\in\X^n: \abs{f(x)-\mu f}>t}}.
\eeq

The 
%normalized 
(unnormalized)
Hamming metric on $\X^n$ is defined by
\beqn
\label{eq:ham}
d(x,y) &=& 
%\oo n
\sum_{i=1}^n \pred{x_i\neq y_i},
\qquad x,y\in\X^n,
\eeqn
where the indicator variable
$\pred{\cdot}$ 
assigns 0-1 truth values 
to the predicate in 
$\set{\cdot}$. 
%The {\em normalized} Hamming metric is $d(x,y)/n$.

The Lipschitz constant of a function, with respect to some metric $d$, is defined by
\beq
\Lip{f} &=& \sup_{x\neq y\in\X^n} \frac{\abs{f(x)-f(y)}}{d(x,y)}.
\eeq

Random variables are capitalized ($X$), specified sequences
are written in lowercase ($x\in\X^n$), the shorthand
$\sseq{X}{i}{j}\defeq
\seq{X}{i}{j}
$ is used for all sequences, and
sequence concatenation is denoted multiplicatively:
$\sseq{x}{i}{j}\sseq{x}{j+1}{k}=\sseq{x}{i}{k}$.
Sums will range over the entire space of the summation variable;
thus
$\ds\sum_{\sseq{x}{i}{j}}f(\sseq{x}{i}{j})$ stands for
$$\ds\sum_{\sseq{x}{i}{j}\in\X^{j-i+1}}f(\sseq{x}{i}{j}).$$ 
By convention, when $i>j$, we define
$$ \sum_{\sseq{x}{i}{j}}f(\sseq{x}{i}{j})
\equiv 
f(\eps)$$
where $\eps$ is the null sequence. 
Products of spaces and measures are denoted by $\tp$.

The {\em total variation} norm of a signed measure $\nu$ on $\X^n$ (i.e., vector 
%in 
$\nu\in\R^{\X^n}$) is defined by
\beq
\TV{\nu} = {\ts\oo2}\nrm{\nu}_1
= {\ts\oo2}\sum_{x\in\X^n} \abs{\nu(x)}
\eeq
(the factor of $1/2$ is not entirely standard).
%note the non-standard factor of $1/2$. 
For readability, we will drop the subscript TV from the norm;
thus everywhere in the sequel, $\nrm{\cdot}$ will mean $\TV{\cdot}$.

A signed measure $\nu$ on a set $\calX$ is called {\em balanced} if $\nu(\calX)=0$. 
Departing from standard convention, our stochastic matrices will be column- (as opposed to row-) stochastic. 
%Definitions for tensors are given in the 
%text
%sequel.
We will use $\abs{\cdot}$ to denote set cardinalities,
and write
$\nn{n}$ for the set $\{1,\ldots,n\}$.

\renewcommand{\TV}[1]{\nrm{#1}}

\section{Contraction and tensorization}
\label{sec:contr}
Our method for bounding $\eta$-mixing coefficients rests on the following simple result:
\belen
%[Markov, 1906 \cite{mar1906}]
\label{lem:markcontr}
Let
%$P$ an $m\times m$ column-stochastic matrix\\
$P:\R^\X\to\R^\X$ be a Markov operator:
$$ (P\nu)(x) = \ds\sum_{y\in\X}P(x\gn y)\nu(y), $$
where $P(x\gn y)\geq0$ and $\sum_{x\in\X} P(x\gn y)=1$.
Define the 
%Doeblin 
contraction
coefficient of $P$ as above:
\beq
\tha 
&=&\max_{y,y'\in\X}\TV{P(\cdot\gn y)-P(\cdot\gn y')}.
\eeq
Then
\beq
\TV{P\nu} \;\leq\; \tha\TV{\nu}
\eeq
for any balanced signed measure
$\nu$ on $\X$
(i.e., $\nu\in\R^\X$ with $\sum_{x\in\X} \nu(x)= 0$).
\enlen
This result is sometimes credited to Dobrushin \cite{dobrushin56}; 
the quantity $\tha$ has been referred to in the literature, alternatively, as the {\em Doeblin contraction}
or {\em Dobrushin ergodicity} coefficient.
However, the 
%inequality 
observation
apparently goes as far back as Markov himself 
\cite{mar1906}
%As indicated, this has been known since Markov himself; 
(see \cite{kontram06} for a proof),
so it seems appropriate to refer to the result above as the
{\em Markov contraction lemma}.

Another important property of the total variation norm is that it tensorizes, 
%as follows:
in the following way:
\belen
\label{lem:tvtens}
Consider two
finite sets $\calX,\calY$, with probability measures $p,p'$ on $\calX$ and
$q,q'$ on $\calY$.
Then 
\beq
\TV{p\tp q - p'\tp q'} &\leq &
\TV{p-p'} + \TV{q-q'}
- \TV{p-p'}\TV{q-q'}.
\eeq
\enlen
This 
%result 
fact
seems to be folklore knowledge; we were not able to locate it in published literature. A non-coupling proof is given in \cite{kont07-thesis}, but we will give
the simpler coupling proof here.
\bepf[Proof of Lemma \ref{lem:tvtens}]
Recall that $$\TV{r-r'}=\inf\pr{Z\neq Z'}$$
where the infimum is over all {\em couplings} of $r$ and $r'$ 
--- i.e., joint distributions on $(Z,Z')$ with respective 
marginals $r$ and $r'$. 
A coupling achieving the infimum is called {\em optimal}. 
Let $X,X',Y,Y'$ be random variables with distributions    
    $p,p',q,q'$, respectively.
Define $\pi$ be an optimal coupling of $p$ and $p'$,
and define similarly $\pi'$ for $p',q'$. 
Notice that
$\pi\tp\pi'$ 
is a (not necessarily optimal) coupling 
of
$p\tp q$ and $p'\tp q'$. 
Then
\beq
\TV{p\tp q-q'\tp p'} &=& \inf\pr{ (X,Y)\neq(X',Y')}\\
&=& \inf\sqprn{\pr{X\neq X'}+\pr{Y\neq Y'}-\pr{X\neq X',Y\neq Y'}}\\
&\le& 
\pi\set{X\neq X'} 
+ 
\pi'\set{Y\neq Y'}
-\pi\set{X\neq X'}\pi'\set{Y\neq Y'} \\
&=&  \TV{p-p'} + \TV{q-q'}
- \TV{p-p'}\TV{q-q'}.
\eeq
\enpf

\section{Markov chains}
\subsection{Directed}
\label{sec:mc}
Technically, this section might be considered superfluous, since this result has already appeared in 
\cite{kontram06}, and is strictly generalized in 
later sections. However, we find it instructive to work out the simple Markov case
as it provides the cleanest illustration of our technique. 

Let $\mu$ be an inhomogeneous Markov measure on $\X^n$, induced by the
kernels $p_0$ and $p_i(\cdot\gn\cdot)$, $1\leq i<n$. Thus,
\beq
\mu(x) &=&  p_0(x_1)\prod_{i=1}^{n-1} p_i(x_{i+1}\gn x_i).
\eeq

The $i$\th contraction coefficient, $\tha_i$ is defined as in (\ref{eq:thadef}). As stated in the Introduction, Markov
contraction 
provides an estimate on
%may be used
%turns out 
%to 
%control 
estimate
$\eta$-mixing:
\beth
\label{thm:marketa}
\beq
\etab_{ij} &\leq& \tha_i\tha_{i+1}\ldots\tha_{j-1}.
\eeq
\enth
\bepf
Fix $1\leq i<j\leq n$ and $\sseq{y}{1}{i-1}\in\X^{i-1}$,
$w_i,w_i'\in\X$. Then
\beq 
\eta_{ij}(y,w,w') 
&=& 
{\ts\oo2}
\sum_{\sseq{x}{j}{n}
}\abs{
\mu(\sseq{x}{j}{n}\gn \sseq{y}{1}{i-1}w_i)
-    
\mu(\sseq{x}{j}{n}\gn \sseq{y}{1}{i-1}w_i')
}\\
&=& {\ts\oo2}
\sum_{\sseq{x}{j}{n}}
\pi(\sseq{x}{j}{n})\abs{
%\zeta(\sseq{x}{j}{n},\sseq{y}{1}{i},w_i,w_i')
\zeta(
%\sseq{x}{j}{n}
x_j
)
}
\eeq
where
\beq
\pi(\sseq{u}{k}{l}) &=& \prod_{t=k}^{l-1} p_{t}(u_{t+1}\gn u_{t})
\eeq
and
\beqn
\label{eq:zdef}
%\zeta(\sseq{x}{j}{n},\sseq{y}{1}{i-1},w_i,w_i')
\zeta(
%\sseq{x}{j}{n}
x_j
)
&=&
\left\{
\begin{array}{ll}
\ds\sum_{\sseq{z}{i+1}{j-1}}
p_{j-1}(x_{j}\gn z_{j-1})
\pi(\sseq{z}{i+1}{j-1})
\paren{
p_i(z_{i+1}\gn w_i)
-
p_i(z_{i+1}\gn w_i')
}
,& j-i>1 \\\\
p_i(x_j\gn w_i)
-
p_i(x_j\gn w_i')
,& j-i=1.
%\quad \mbox{otherwise}.
\end{array}
\right.
\eeqn
%Define $\vec x\in\R^\X$ by
%$$\vec x_v 
%= p_{i-1}(v \gn y_{i-1})
%.$$ 

Define $\vec h\in\R^\X$ by 
$\vec h_v = p_i(v\gn w_i) - p_i(v\gn w_i')$
and $P\supr k\in\R^{\X\times\X}$ by
$P\supr k_{u,v} = p_k(u\gn v)$.
Likewise, define
$\vec z\in\R^\X$ by 
$\vec z_v = 
%\zeta([v\,\sseq{x}{j+1}{n}],\sseq{y}{1}{i-1},w_i,w_i')
%\zeta([v\,\sseq{x}{j+1}{n}])
\zeta(v)
$.
It follows that
\beq
\vec z = 
P\supr{j-1}
P\supr{j-2}
%\ldots
\cdots
P\supr{i+2}
P\supr{i+1}
\vec h.
\eeq
Therefore,
\beq
\eta_{ij}(y,w,w') 
&=&
{\ts\oo2}
\sum_{\sseq{x}{j}{n}}
\pi(\sseq{x}{j}{n})\abs{\vec z_{x_j}}\\
&=&
{\ts\oo2}
\sum_{x_j}
\abs{\vec z_{x_j}}
\sum_{\sseq{x}{j+1}{n}}
\pi(\sseq{x}{j}{n})\\
&=&
{\ts\oo2}
\sum_{x_j}
\abs{\vec z_{x_j}}
= \TV{\vec z}.
\eeq
The claim follows by (repeated applications of) the Markov contraction lemma.
\enpf

The reader may wish to compare this proof with Samson's \cite{samson00}.
\begin{rem}
Aside from simplicity, another advantage of the Markov contraction method is its precision.
Namely, in the proof above, note that equality is maintained until the very end, where the Markov contraction
lemma is invoked to bound $\TV{\vec z}$. This means that more delicate (for example, spectral) 
estimates on $\TV{
P\supr{j-1}
P\supr{j-2}
%\ldots
\cdots
P\supr{i+2}
P\supr{i+1}
\vec h}$ translate directly into tighter concentration bounds.
\end{rem}

\subsection{Undirected}
\renewcommand{\a}{\alpha}
\renewcommand{\b}{\beta}

In this section we analyze Markov chains under a different parametrization, in an ``undirected graphical
model'' setting \cite{lauritzen96}.
For any graph $G=(V,E)$, where $|V|=n$ and the maximal cliques have size 2 (i.e., are
edges), we can define a measure on $\X^V=\X^n$ as follows
\beqn
\label{eq:und-mudef}
\mu(x)
= \frac{\prod_{(i,j)\in E} \psi_{ij}(x_i,x_j)}
{\sum_{x'\in\X^n}\prod_{(i,j)\in E} \psi_{ij}(x'_i,x'_j)}
%\equiv
%\frac{\prod_{(i,j)\in E} \psi_{ij}(x_i,x_j)}
%{Z_G}
,\qquad x\in\X^n
\eeqn
%There is no loss of generality in taking $\psi_{ij}(x,y)\geq1$, which
%we do henceforth.
for some 
for some nonnegative ``potential functions'' $\psi_{ij}$.
%$\psi_{ij}\geq0$.

Consider the very simple case of chain graphs; any such measure is a
Markov measure on $\X^n$.
% (as before, $|V|=n$). 
We can relate the induced Markov transition
kernel $p_i(\cdot\gn\cdot)$ to the random field measure $\mu$ as follows:
\beq
p_i(x\gn y) &=&
\frac{
\sum_{\sseq{v}{1}{i-1}}\sum_{\sseq{z}{i+2}{n}} \mu(v y x z)
}{
\sum_{x'\in\X}\sum_{\sseq{{v'}}{1}{i-1}}\sum_{\sseq{{z'}}{i+2}{n}}
 \mu(v' y x' z')
},\qquad
x,y\in\X
.
\eeq
%(Note: expressions such as $\sum_{\sseq{u}{i}{j}}$ should be
%interpreted as $\sum_{\sseq{v}{i}{j}\in\X^{j-i+1}}$; in general, the
%subscripts and superscripts on a word indicate its length and
%therefore the set being summed over.)

Our goal is to bound the $i$\th contraction coefficient $\theta_i$ of the Markov
chain\hide{
:
\beq
\theta_i &=& \max_{y,y'\in\X} \oo2
\sum_{x\in\X}\abs{p_i(x\gn
  y)-p_i(x\gn y')}.
\eeq
%and the quantity $\rho_i$:
%\beq
%\rho_i &=& 
%\max_{x,y\in\X}\psi_{i,i+1}(x,y).
%\eeq
}
in terms of $\psi_{ij}$.
We claim 
%that there is 
a simple relationship between $\theta_i$ and $\psi_{ij}$:
\bethn
\label{thm:undmark}
\beqn
\label{eq:thrho}
\theta_i &\leq& \frac{R_i-r_i}{R_i+r_i},
\qquad 1\leq i<n
\eeqn
where
$$ R_i=\max_{x,y\in\X}\psi_{i,i+1}(x,y) $$
and
$$ r_i=\min_{x,y\in\X}\psi_{i,i+1}(x,y). $$
\enthn

%First we prove a simple lemma:
We will need the following lemma, which may be of independent interest:
\belen
\label{lem:auv}
For $n\in\N$
and $0\le r\le R$,
consider the vectors 
$\a\in[0,\infty)^n$
and
$f,g\in[r,R]^n$.
Then
\beq
{\oo2}\sum_{i=1}^{n}\abs{
\frac{\a_if_i}{\sum_{j=1}^{n} \a_jf_j}
-
\frac{\a_ig_i}{\sum_{j=1}^{n} \a_jg_j}
} &\leq& \frac{R-r}{R+r}.
\eeq
\enlen
\begin{rem*}
This lemma was proved together with Roi Weiss.
\end{rem*}
\bepf
Assume for now $r>0$ (the general case will follow by continuity).
Since 
$$\sum_{i=1}^n|x_i-y_i|=\max_{\beta\in\set{-1,1}^n}\sum_{i=1}^n\beta_i(x_i-y_i),$$
proving our claim is equivalent to showing that
\beqn
\label{eq:abxy}
\sum_{i}\b_i\paren{
\frac{\a_if_i}{\sum_{j} \a_jf_j}
-
\frac{\a_ig_i}{\sum_{j} a_jg_j}
} &\leq& 2\frac{R-r}{R+r}
\eeqn
holds for all $f,g\in[r,R]^n$,
$\a\in[0,\infty)^n$,
and $\b\in\set{-1,1}^n$. 
For fixed $\a,\b$, define the function $F_{\a,\b}:[r,R]^{2n}\to\R$ by putting 
$F_{\a,\b}(f,g)$ equal to the left-hand side of (\ref{eq:abxy}).
If $F_{\a,\b}$ achieves an extremum somewhere on $[r,R]^{2n}$, its gradient must vanish there.
Let us compute this gradient:
\beq
\ddel{F_{\a,\b}}{f_i} &=& \frac{ \a_i\sum_{j\neq i}\a_j(\b_i-\b_j)f_j}{\paren{\sum_{k=1}^n\a_i f_i}^2} \\
\ddel{F_{\a,\b}}{g_i} &=& \frac{ \a_i\sum_{j\neq i}\a_j(\b_i-\b_j)g_j}{\paren{\sum_{k=1}^n\a_i g_i}^2}.
\eeq
Solving for $\grad{F_{\a,\b}}\equiv0$, we get that the latter holds whenever
\beq
f_a = \frac{\sum_{j\notin\set{a,b}} \a_j(\b_b-\b_j)f_j}{\a_a(\b_a-\b_b)},
\qquad
f_b = \frac{\sum_{j\notin\set{a,b}} \a_j(\b_j-\b_a)f_j}{\a_b(\b_a-\b_b)},\\
g_c = \frac{\sum_{j\notin\set{c,d}} \a_j(\b_d-\b_j)g_j}{\a_c(\b_c-\b_d)},
\qquad
g_d = \frac{\sum_{j\notin\set{c,d}} \a_j(\b_j-\b_c)g_j}{\a_d(\b_c-\b_d)}
\eeq
for some $a\neq b$ and $c\neq d$.
Since the expressions above are undefined for $\b_a=\b_b$ or $\b_c=\b_d$,
we may assume that neither of these holds. We claim that
$f_a\le0$ 
for $a\neq b$ and $\b_a\neq\b_b$; this is easily verified by
substituting $\b_a=1,\b_b=-1$ and $\b_a=-1,\b_b=1$. (A similar observation
holds for $g_c$.) We conclude that $\grad F_{\a,\b}$ does not vanish 
anywhere on 
%the interior of 
$[r,R]^{2n}$, and therefore the function must achieve its
extreme values on the boundary of this region.

In light of the above, it suffices to consider $f,g\in\set{r,R}^n$. Consider the expression
\beqn
\label{eq:ABuv}
\calL =
\sum_{i}\b_i\paren{
\frac{\a_if_i}{A}
-
\frac{\a_ig_i}{B}
},
\eeqn
where $f,g\in\set{r,R}^n$, $\a\in[0,\infty)^n$, $\b\in\set{-1,1}^n$, and
$A=\sum_j \a_jf_j$, $B=\sum_j \a_jg_j$. Keeping $A$ and $B$ constant, we seek an $\a$ that maximizes 
%(\ref{eq:ABuv}).
$\calL$.
The latter is a linealy constrained linear program, and thus its maximal value(s) are
attaind at the extreme points of the feasible region $E=E(A,B,f,g)\subset\R^n$, given by
$$ E=\set{\a\in[0,\infty)^n : \iprod{\a,f}=A, \iprod{\a,g}=B}.$$
Recalling that 
$z\in E$ is
an extreme point 
if and only if $z$ cannot be expressed as $z=\lambda x+(1-\lambda)y$ for $x,y\in E$ and $0<\lambda<1$,
it is straightforward to verify that the extreme points of $E$ are necessarily
of the form
$\a_i=0$ for  $i\notin\set{a,b}$ for some $a\neq b\in[n]$.
%\beq
%A = \a_a f_a + \a_b f_b + \sum_{i\notin \{a,b\}} \a_i f_i\\
%B = \a_a g_a + \a_b g_b + \sum_{i\notin \{a,b\}} \a_i g_i,
%\eeq
%for all $a\neq b$ pairs and so the extermal points are with $\a_i = 0$ for $i\notin \{a,b\}$.
Furthermore, since $\calL$ is homogeneous in $\a$, we may rescale it so that
 $\a_b=1$.
In this case we have
\beq
\a_a=\frac{A-f_b}{f_a}
\eeq
and
\beq
B=\a_a g_a + g_b = \frac{A-f_b}{f_a} g_a + g_b.
\eeq
Substituting this into
(\ref{eq:ABuv}), after some algebra we get
\beq
%\label{eq:AB2}
\calL=\frac{(\b_a - \b_b) (A-f_b)(g_b f_a - g_a f_b)}
	  {A(A g_a + g_b f_a - g_a f_b)}.
\eeq
Maximizing $\calL$ with respect to $A$ yields
%Since we are looking for a maximum let us maximize the last display w.r.t $A$. Differentiating and equating to zero we get a quadratic equation with the solution
\beq
A^* = f_b + \sqrt{\frac{f_a f_b g_b}{g_a}},
\eeq
with a maximal value of
%Putting this solution back to (\ref{eq:AB2}) we get the expression
\beqn
\label{eq:LA*}
\calL(A^*) = 
(\b_a - \b_b) \frac
	{\left(
				g_b f_a - g_a f_b
	 \right)
	 \sqrt{\frac{f_a f_b g_b}{g_a}}}
	{
	 \left(
			f_b + \sqrt{\frac{f_a f_b g_b}{g_a}}
	 \right)
	 \left(
			g_b f_a + g_a \sqrt{\frac{f_a f_b g_b}{g_a}}
	 \right) 
	}.
\eeqn
The nontrivial values for $\b_a,\b_b$ are $1$ and $-1$, respectively, and it remains
%Setting $\b_a = 1$ and $\b_b = -1$ we gauranty the maximum and all we have to do now is 
%to check what values of 
to choose the
$f_a,f_b,g_a,g_b \in \{r,R\}$ 
so as to maximize the above display.
%make the above expression largest. 
%Setting 
It is easily seen that the choice
$f_a = g_b = R$ and $f_b = g_a = r$ 
%achieve the task with
is optimal, which yields the value
\beq
2\frac{(R^2 - r^2) R }{R(r + R)^2} = 
2\frac{R-r}{R+r},
\eeq
as claimed.
\enpf
\begin{rem}
In cases where the estimate $(R-r)/(R+r)$ is too crude, 
it might be possible to salvage a
useful bound from (\ref{eq:LA*}).
\end{rem}

\bepf[Proof of Theorem \ref{thm:undmark}]
Let us define the shorthand notation:
\beq
\pi(\sseq{u}{k}{l}) &=& \prod_{t=k}^{l-1} 
\psi_{t,t+1}(u_t,u_{t+1})
\eeq
Then we expand
\beq
p_i(x\gn y)
&=&
\frac{
\sum_{\sseq{v}{1}{i-1}}\sum_{\sseq{z}{i+2}{n}} 
\pi(\sseq{v}{1}{i-2})\psi_{i-1,i}(v_{i-1},y)\psi_{i,i+1}(y,x)\psi_{i+1,i+2}(x,z_{i+2})\pi(\sseq{z}{i+2}{n})
}{
\sum_{x'\in\X} \sum_{\sseq{{v'}}{1}{i-1}}\sum_{\sseq{{z'}}{i+2}{n}}
\pi(\sseq{{v'}}{1}{i-2})\psi_{i-1,i}(v'_{i-1},y)\psi_{i,i+1}(y,x')\psi_{i+1,i+2}(x',z'_{i+2})\pi(\sseq{{z'}}{i+2}{n})
}\\
&=&
\frac{
\psi_{i,i+1}(y,x)a_{yx}
}{
\sum_{x'\in\X} \psi_{i,i+1}(y,x')a_{yx'}
}
\eeq
where
\beq
a_{yx} &=&
\sum_{\sseq{v}{1}{i-1}}\sum_{\sseq{z}{i+2}{n}} 
\pi(\sseq{v}{1}{i-2})\psi_{i-1,i}(v_{i-1},y)\psi_{i+1,i+2}(x,z_{i+2})\pi(\sseq{z}{i+2}{n})
\eeq
(we take the natural convention that 
$\psi_{i,j}(\cdot\gn\cdot)=1$ whenever $(i,j)\notin E$).

Fix $y,y'\in\X$. Define the quantities, for each $x\in\X$:
\beq
f_x &=& \psi_{i,i+1}(y,x) \\
g_x &=& \psi_{i,i+1}(y',x) \\
\a_x &=& a_{yx} \\
\a'_x &=& a_{y'x}.
\eeq
Then
\beq
\sum_{x\in\X}\abs{p_i(x\gn
  y)-p_i(x\gn y')}
&=&
\sum_{x\in\X}\abs{
\frac{f_x\a_x}{\sum_{x'\in\X} f_{x'}\a_{x'}}
-
\frac{g_x\a'_x}{\sum_{x'\in\X} g_{x'}\a'_{x'}}
}\\
&=&
\sum_{x\in\X}\abs{
\frac{f_x\a_x}{\sum_{x'\in\X} f_{x'}\a_{x'}}
-
\frac{g_x\a_x}{\sum_{x'\in\X} g_{x'}\a_{x'}}
};
\eeq
the last equality follows since $\a'_x=c\a_x$, where
$c=\frac{\psi_{i-1,i}(v_{i-1},y')}{\psi_{i-1,i}(v_{i-1},y)}$.
Now Lemma \ref{lem:auv} can be applied to establish the claim.
\enpf
%Since all the inequalities invoked are tight, so is the bound in
%Theorem \ref{thm:undmark}.

\hide{
%Let us compare 
Theorem \ref{thm:undmark} 
compares 
%rather 
favorably to
some of the results in the literature. 
}
We observe that in general, mixing bounds on directed Markov chains are more informative
than on undirected ones.
In the language of random fields,
the measure $\mu$ defined in (\ref{eq:und-mudef}) is a finite volume Gibbs measure with a pair potential.
Let us recall the Dobrushin uniqueness condition and its role in concentration inequalities. For $j\in V$, 
define 
$\omit{V}{j}
:=V\setminus\set{j}$
and
%let 
$\sim_j$ to be the following equivalence relation on $\Omega^V$: $x\sim_j y$ if $x_k=y_k$ for all $k\in\omit{V}{j}$.
Also, define the operator $\omit{(\cdot)}{i}:\Omega^V\to \Omega^{
\omit{V}{i}
}$ as the obvious coordinate projection.
For $i\in V$ and $x\in\Omega^V$, define $\mu_i(\cdot\gn \omit{x}{j}
)$ to be the distribution on $X_i$ conditioned on 
the other $\set{X_j}$ being equal to $\omit{x}{j}$. 
Define the $n\times n$ {\em Dobrushin interdependence matrix} $D=(d_{ij})$:
\beqn
\label{eq:dobr-d}
d_{ij} = \max_{x \sim_j y \in \Omega^V}\TV{ \mu_i(\cdot\gn \omit{x}{j}) - \mu_i(\cdot\gn \omit{y}{j}) }.
\eeqn
The Dobrushin uniqueness condition requires that $\nrm{D}_\infty<1$ (all norms on 
discussed here
$D$ are $\ell_p\to\ell_p$ operator norms).
Let $\mu$ be any probability measure on $\Omega^V$ and for simplicity, take $f:\Omega^V\to\R$ to 
satisfy $\sqrt n\Lip{f}\le1$
with respect to the Hamming metric.
Then a result of K{\"u}lske \cite[Theorem 1]{kulske03} states that
\beqn
\label{eq:kulske}
\mu\set{ |f-\mu f| \ge t} \le 2\exp\paren{-\frac{t^2}{2} (1-\dsnrm{D}_\infty)(1-\dsnrm{D\trn}_\infty)},
\qquad t\ge0,
\eeqn
provided that $\dsnrm{D}_\infty,\dsnrm{D\trn}_\infty<1$,
and Chatterjee in \cite[Theorem 4.3]{chatterjee05} proves that
\beqn
\label{eq:chatter}
\mu\set{ |f-\mu f| \ge t} \le 2\exp\paren{-(1-\nrm{D}_2)t^2}
\qquad t\ge0,
\eeqn
provided that $\nrm{D}_2<1$.

On the other hand, consider a directed
homogeneous
Markov chain $(X_1,\ldots,X_n)$ with measure $\mu$ on $\Omega^n$ induced by a Markov kernel with 
contraction coefficient $\theta<1$. 
In this case, the $\ell_\infty$ norm of the $\eta$-mixing matrix (see \ref{eq:Delta-inf}) is easily bounded:
\beq
\nrm{\Delta}_\infty \le \oo{1-\theta}.
\eeq

An immediate consequence of the definition (\ref{eq:dobr-d}) is that $d_{n,n-1}=\theta$, which implies that
\beqn
\label{eq:thaD}
\theta \le 
\min\set{
\dsnrm{D}_\infty,\dsnrm{D\trn}_\infty,\nrm{D}_2
}.
\eeqn
Thus, for directed homogeneous Markov chains,
(\ref{eq:kontram}) is uniformly superior to (\ref{eq:kulske}).
The inequality (\ref{eq:thaD}) does not imply that
(\ref{eq:kontram}) is uniformly superior to (\ref{eq:chatter}) 
--- and indeed, this is not the case.
However, (\ref{eq:thaD}) does imply that whenever 
%Chatterjee's bound 
(\ref{eq:chatter}) 
is nontrivial (i.e., $\nrm{D}_2<1$),
the contraction
%our 
bound (\ref{eq:kontram}) is also nontrivial (i.e., $\theta<1$).
This implication does not hold in the reverse, as we show by example.
Take $\Omega=\set{0,1}$ and consider the Markov kernel defined by
$P(0\gn0)=a$, $P(0\gn1)=b$,
$P(1\gn0)=1-a$, $P(1\gn1)=1-b$, for $a,b\in(0,1)$. A straightforward calculation 
shows that 
for $a\le1/2$ and
as $b\searrow0$, we have $\theta\to a$ and 
$$\min\set{
\dsnrm{D}_\infty,\dsnrm{D\trn}_\infty
}\to 2/(1+a),\qquad
\nrm{D}_2>1.5/(1+a),$$
rendering the corresponding 
Dobrushin matrix-based
estimates uninformative
even for very small $\theta$.

%$d_{1,2}=d_{n,n-1}=\theta=

\section{Markov tree processes}
\label{sec:mt}
\subsection{Preliminaries}
\newcommand{\lev}{\operatorname{lev}}
\newcommand{\rents}{\operatorname{parents}}
\newcommand{\kids}{\operatorname{children}}
\newcommand{\depth}{\operatorname{dep}}
\newcommand{\width}{\operatorname{wid}}
\newcommand{\dpath}{\vec \pi}
\newcommand{\rsub}[1]{_{#1}}  % real subscript
\newcommand{\bsub}[1]{[#1]}   % bracket subscript
\newcommand{\inex}{\alpha}

%The material in this section is taken almost entirely from
%\cite{kont06-tree}.
We begin by defining
%Let us first define 
some notation specific to this section. A
collection of variables may be indexed by 
%Another way to index collections of variables is by 
subset: if
$x\in\X^V$ and $I\subseteq V$ with
$I=\{i_1, i_2,\ldots, i_m\}$, then we write 
$x_I\equiv x[I]\defeq \{x_{i_1},x_{i_2},\ldots,x_{i_m}\}$;
we will write $x_I$ and $x[I]$ interchangeably, as dictated by
convenience.
To avoid cumbersome subscripts, we will also occasionally use the
bracket notation for vector components. 
Thus, 
if
%if $I\subset\NN$
%$I\subset V$ and 
$\vec u\in \R^{\X^I}$, then
$$ \vec u_{x_I}\equiv \vec u_{x[I]} \equiv \vec u[x_I] \equiv \vec u[x[I]]
\defeq \vec u_{(x_{i_1},x_{i_2},\ldots,x_{i_m})} \in \R$$
for each $x[I]\in\X^I$.
A similar bracket notation will apply for matrices.
If $A$ is a matrix then $A_{*,j}=A[*,j]$ will denote its $j$\th column.
Probabilities are denoted by $\P$ in this section.

If $G=(V,E)$ is a graph,
we will 
frequently
%sometimes 
abuse notation and write $u\in G$ instead of $u\in
V$, blurring the distinction between a graph and its vertex set. This
notation will carry over to set-theoretic operations ($G=G_1\cap G_2$)
and indexing of variables (e.g., $X_G$).
%To
%specify the edge set of a graph, we will write $E(G)$.

\subsection{Graph theory}
\label{sec:graph}
%We shall need some
%graph-theoretic preliminaries.
Consider a directed acyclic graph 
%(DAG) 
$G=(V,E)$, and define a
partial order 
$\prec_G$
on $G$ by the transitive closure of the relation
$$u\prec_G v 
%for each directed edge 
\qquad\text{if}\qquad
(u,v)\in E.$$
We define the {\em parents} and {\em children} of $v\in V$ in the
natural way:
$$ \rents(v) = \{u\in V : (u,v)\in E\} $$
and
$$ \kids(v) = \{w\in V : (v,w)\in E\}. $$

If $G$ is connected and each $v\in V$ has at most one parent, $G$ is called a
{\em(directed) tree}. In a
tree, whenever $u\prec_G v$ there is a unique directed path from $u$ to
$v$.
%; denote its length (i.e., the number of edges it traverses) by
%$\dpath(u,v)$.
A tree $T$ always has a unique minimal (with respect to
%under 
$\prec_T$)
element $r_0\in V$, called its {\em root}. Thus, for every $v\in V$
there is a unique directed path 
$r_0 \prec_T r_1 \prec_T \ldots \prec_T r_d = v$; define
the {\em depth} of $v$,
$\depth_T(v)=d$, to be the length (i.e., number of edges) of this path.
Note that $\depth_T(r_0)=0$. We define the depth of the tree by
$\depth(T)=\sup_{v\in T}\depth_T(v)$.

For $d=0,1,\ldots$ define
the $d$\th {\em level} of the tree $T$ by
$$ \lev_T(d) = \{v \in V : \depth_T(v)=d\};$$
note that the levels induce a disjoint partition on $V$:
$$ V = \bigcup_{d=1}^{\depth(T)} \lev_T(d).$$
We define the {\em width}\footnote{
%Note that t
This definition is nonstandard.
} of a tree as the greatest number of nodes in
any level:
\beqn
\label{eq:wid}
 \width(T) = \sup_{1\leq d\leq\depth(T)} \abs{\lev_T(d)}.
\eeqn

We will consistently take $|V|=n$ for finite $V$.
An ordering 
%$J:V\to\nn{n}$
$J:V\to\N$
of the nodes is 
%called 
said to be
{\em breadth-first} if 
\beqn
\label{eq:depj}
 \depth_T(u) < \depth_T(v) \Longrightarrow J(u) < J(v) .
\eeqn
Since every 
finite
directed tree
$T=(V,E)$ has some breadth-first ordering,\footnote{
One can easily construct a breadth-first ordering on a given tree by
ordering the nodes arbitrarily within each level and listing the
levels in ascending order: $\lev_T(1),\lev_T(2),\ldots$.
}
we will henceforth blur
the distinction between $v\in V$ and $J(v)$, 
simply taking 
%$V=\{1,2,\ldots,n\}$
$V=\nn{n}$
(or $V=\N$)
and assuming that 
%(\ref{eq:depj})
$\depth_T(u) < \depth_T(v) \Rightarrow u < v$ 
holds. This
will allow us to write $\X^V$ simply as $\X^n$
for any set $\X$.

Note that we have two orders on $V$: the partial order $\prec_T$,
induced by the tree 
%topology, 
edges,
and the total order $<$, given by the
breadth-first enumeration. Observe that $i\prec_T j$ implies $i<j$ but
not 
vice versa.
%the other way around. 
%Another observation is that this setup
%forces $i=1$ to be the root node.

If $T=(V,E)$ is a tree and $u\in V$, we define the {\em subtree} 
induced by $u$,
$T_u=(V_u,E_u)$
by $V_u = \{v\in V: u\preceq_T v\}$, 
$E_u = \{(v,w) \in E: v,w\in V_u\}$.

\subsection{Markov tree measure}
\label{sec:MTmeas}
If $\X$ is a finite set,
%and $\calF=2^\X$, 
a
{\em Markov tree measure} $\mu$ is defined on
%$(\X^n,\calF^n)$ 
$\X^n$ 
by a tree $T=(V,E)$
and transition kernels $p_0$, 
$\tlprn{p_{ij}: (i,j)\in E}$.
Continuing our convention above,
%in \S\ref{sec:graphtheo}, 
we have a
breadth-first order $<$ and the total order $\prec_T$ on $V$, and take
$V=\{1,\ldots,n\}$. Together, the 
%topology 
edges
of $T$ and the transition
kernels determine the 
%measure 
distribution
$\mu$ on $\X^n$:
\beqn
\label{eq:mtmeas}
\mu(x) = p_0(x_1)\prod_{(i,j)\in E} p_{ij}(x_j\gn x_i),
\qquad x\in\X^n.
\eeqn
A measure
%Furthermore, we refer to any measure 
%$\nu$ 
on 
$\X^n$ 
satisfying
%that can be
%written in the form of 
(\ref{eq:mtmeas}) for some $T$ and $\{p_{ij}\}$
is said to be {\em compatible} with tree $T$;
a measure is
%as 
a Markov tree measure if it is compatible with some tree.

Suppose $\X$ is a finite set and $(X_i)_{i\in\N}$, $X_i\in\X$ is a
random process 
defined on 
%$(\X^\N,\calF^\N,\P)$. 
$(\X^\N,\P)$. 
If for each $n>0$
there is a tree $T\supr n=(\nn{n},E\supr n)$ 
and a Markov tree measure $\mu_n$
compatible with $T\supr n$
such that
for all $x\in\X^n$ we have
$$ \pr{\sseq{X}{1}{n} = x} = \mu_n(x)$$
%for some Markov tree measure $\mu_n$ on 
%$(\X^n,\calF^n)$,
%$\X^n$
then we call $X$ a {\em Markov tree process}. 
The trees $\{T\supr n\}$ are easily seen to be consistent in the sense
that $T\supr n$ is an induced subgraph of $T\supr{n+1}$.
So corresponding 
to 
%this
any Markov tree process is the unique infinite tree
%process is the infinite tree 
$T=(\N,E)$. 
%The latter 
The uniqueness of $T$ 
is easy to see,
%This $T$ is unique,
since for
$v>1$, the parent of $v$ is the 
%unique 
smallest
$u\in \N$ such that
$$ \pr{X_v = x_v \gn \sseq{X}{1}{u}\;=\;\sseq{x}{1}{u}} = \pr{X_v = x_v \gn X_u=x_u};$$
thus $\P$ determines the 
%topology 
edges
of $T$. 

It is straightforward to
verify 
%(via (\ref{eq:mtmeas}))
that 
a Markov tree process $\{X_v\}_{v\in T}$
compatible with tree $T$
has the following {\em Markov property}: if $v$ and $v'$ are children
of $u$ in $T$, then
$$ \pr{X\rsub{T_v}=x,X\rsub{T_{v'}}=x'\gn X_{u}=y}
= \pr{X\rsub{T_v}=x\gn X_{u}=y}
\pr{X\rsub{T_{v'}}=x'\gn X_{u}=y}.$$
In other words, the subtrees induced by the children are conditionally
independent given the parent; this follows directly from the
definition of the Markov tree measure in (\ref{eq:mtmeas}).

\subsection{Statement of result}
\bethn
\label{thm:mtmain}
Let $\X$ be a finite set and
let $(X_i)_{1\leq i\leq n}$, $X_i\in\X$  be a Markov tree process, defined by a
tree $T=(V,E)$ and transition kernels $p_0$, 
$\tlprn{p_{uv}(\cdot\gn\cdot)}_{(u,v)\in E}$.
Define the $(u,v)$- {\em
%(Doeblin-) 
contraction coefficient} $\tha_{uv}$ by
\beqn
\label{eq:uvthadef}
\tha_{uv} &=&
\max_{y,y'\in\X}
\TV{
p_{uv}(\cdot\gn y)-p_{uv}(\cdot\gn y')
}.
\eeqn
Suppose 
$\max_{(u,v)\in E} \tha_{uv}\leq \tha <1 $
for some $\tha$ 
and
$\width(T)\leq L$
\hide{
$$\max
%_{\el:~ \depth(u)=\el,~ 
%1\leq u \leq n
%u\in\nn{n}
%}
\tlprn{
\abs{\lev(d)} 
:
d=\depth(u), u\in\nn{n}
}
\;\leq\; L
$$
for some $L$
(that is, any level of $T$ contains at most $L$ nodes). 
}. 
Then for the
Markov tree process $X$ we have 
\beqn
\label{eq:mainbd}
\etab_{ij} &\leq& \paren{1-(1-\tha)^L}^{\floor{(j-i)/L}}
\eeqn
for $1\leq i<j\leq n$.
\enthn

To cast 
%the result in 
(\ref{eq:mainbd}) in more usable form, we first note that for $k,L\in\N$
%and $k\in\N$, 
with $k\geq L$,
we have
% then
\beqn
\label{eq:flrbd}
 \floor{\frac{k}{L}} \geq \frac{k}{2L-1}
\eeqn
(we omit the elementary number-theoretic proof). Using
(\ref{eq:flrbd}), we have
\beqn
\label{eq:thatil}
\etab_{ij} &\leq& \tilde\tha^{j-i},
\qquad\text{for } j\geq i+L
\eeqn
where
$$ \tilde\tha = (1-(1-\tha)^L)^{
%\oo
1/(
{2L-1}
)
};$$
this implies the 
dimension-free
bound
\beq
\nrm{\Delta}_\infty &\leq& L-1+(1-\tilde\tha)\inv.
\eeq

In the (degenerate) case where
the Markov tree is a chain, we have $L=1$ and therefore
$\tilde\tha=\tha$; thus we recover Theorem \ref{thm:undmark}.
%the Markov chain concentration
%results in \cite{kontram06,marton96,samson00}.
%and the approximations 
%in (\ref{eq:gamapprox},\ref{eq:delapprox})
%become precise inequalities.

\subsection{Proof of 
Theorem \ref{thm:mtmain}
}
\label{sec:proofmain}
The proof of Theorem \ref{thm:mtmain} is combination of elementary graph
theory and tensor algebra. We start with a graph-theoretic lemma:
\belen
\label{lem:j0}
Let $T=(\nn{n},E)$ be a tree and fix $1\leq i<j\leq n$. Suppose 
$(X_i)_{1\leq i\leq n}$ is a Markov tree process 
%for a
%is a 
whose
%measure 
%law
distribution
$\P$
on $\X^n$ 
is
compatible with $T$ 
(in the sense of
%\S\ref{sec:MTmeas}
Section \ref{sec:MTmeas}). 
%(this notion is defined above).
Define the set
$$ 
\hide{
T_{ij} = \{
 k\in T_i :
%T : i\prec_T k, ~
k\geq j\}
}
T_i^j = T_i \cap \{j,j+1,\ldots, n\}
, $$
%that is, $T_{ij}$ is 
consisting of those nodes in the subtree $T_i$ whose
breadth-first numbering does not precede $j$. Then, for $y\in\X^{i-1}$
and $w,w'\in\X$, we have
\beqn
\eta_{ij}(y,w,w') &=&
\left\{
\begin{array}{ll}
0, & \quad T_i^j = \emptyset \\
\eta_{ij_0}(y,w,w'), & \quad \mbox{otherwise},
\end{array}
\right.
\eeqn
where $j_0$ is the minimum (with respect to $<$) element of $T_i^j$.
\enlen
\begin{rem}
This lemma tells us that when computing $\eta_{ij}$ it is sufficient to
restrict our attention to the subtree induced by $i$.
\end{rem}
\bepf
The case $j\in T_i$ implies $j_0=j$ and is trivial; thus we assume
$j\notin T_i$. In this case, the subtrees $T_i$ and $T_j$ are
disjoint. Putting $\bar T_i=T_i\setminus\{i\}$, 
we have by the Markov 
%measure 
property, 
\beq
\pr{X\rsub{\bar T_i}=x\rsub{\bar T_i},X\rsub{T_j}=x\rsub{T_j}\gn 
\sseq{X}{1}{i} = {y}{w}}
&=& 
%\\
%&&\hspace{-1.5cm}
\pr{X\rsub{\bar T_i}=x\rsub{\bar T_i}\gn X_i=w}
  \pr{X\rsub{T_j}=x\rsub{T_j}\gn \sseq{X}{1}{i-1} = y}.
\eeq

%Let 
%$\bar T_i^j = \bar T_i \cap \{j,j+1,\ldots, n\}$;
%$\bar T_i^j = T_i^j\setminus\{i\}$, 
Then
from 
the definition of $\eta_{ij}$
%(\ref{eq:etadef})
% and (\ref{eq:tv}), 
and by
marginalizing out the $X_{T_j}$,
we have
\beq
\eta_{ij}(y,w,w') &=&
{\ts\oo2}
\sumseq{x}{j}{n}
\abs{
\pr{\sseq{X}{j}{n}=\sseq{x}{j}{n}\gn\sseq{X}{1}{i}=yw}
-    
\pr{\sseq{X}{j}{n}=\sseq{x}{j}{n}\gn\sseq{X}{1}{i}=yw'}
}\\
&=&
{\ts\oo2}
\sum_{x\rsub{T_i^j}}
\abs{
\pr{X\rsub{T_i^j}=x\rsub{T_i^j}\gn X_i=w}
-
\pr{X\rsub{T_i^j}=x\rsub{T_i^j}\gn X_i=w'}
}.
\eeq
If $T_i^j=\emptyset$ then obviously $\eta_{ij}=0$; otherwise,
$\eta_{ij}=\eta_{ij_0}$, since $j_0$ is the 
%smallest 
``first''
element of $T_i^j$.
\hide{
This leaves two possible cases:
\bit
\item[(\one)] there is a (smallest) $j_0\in T_i$, $j_0\geq j$
\item[(\two)] for each $k\in T_i$, $k<j$.
\eit
}
\enpf

Next we develop some basic results for tensor norms.
%; recall that
%unless specified otherwise, the
%norm used in this paper is the total variation norm.
If 
%$\vec A\in \R_{M, N}$ is a 
$\vec A$
is an $M\times N$
column-stochastic 
matrix
(i.e., $\vec A_{ij}\geq0$ for $1\leq i\leq M$, $1\leq j\leq N$ and
$\sum_{i=1}^M \vec A_{ij}=1$ for all $1\leq j\leq N$) and $\vec u\in \R^N$
is {\it balanced} in the sense that $\sum_{j=1}^N \vec u_j=0$, we
have, by 
Lemma \ref{lem:markcontr}
\beqn
\label{eq:contr}
\nrm{\vec{Au}} &\leq& \nrm{\vec A}\nrm{\vec u},
\eeqn
where
\beqn
\label{eq:matnorm}
\nrm{\vec A} &=& \max_{1\leq j,j'\leq N} \nrm{\vec A_{*,j}-\vec A_{*,j'}},
\eeqn
and 
$\vec A\rsub{*,j}
\equiv
\vec A\bsub{*,j}
$ 
denotes the $j$\th column of $\vec A$. An immediate
consequence
of (\ref{eq:contr}) is that $\nrm{\cdot}$ 
%is a valid matrix norm:
satisfies
\beqn
\label{eq:AB}
\nrm{\vec{AB}} &\leq& \nrm{\vec A}\nrm{\vec B}
\eeqn
for column-stochastic matrices
$\vec A\in\R^{M\times N}$ and $\vec B\in\R^{N\times P}$.
\begin{rem}
\label{rem:stochnorm}
Note that 
if $\vec A$ is a column-stochastic matrix
then
$\nrm{\vec A}\leq1$, 
and if additionally $\vec u$ is balanced then 
$\vec{Au}$ is also balanced.
\end{rem}

If $\vec u\in\R^M$ and $\vec v\in\R^N$, define their tensor product
$\vec w = \vec v \tp \vec u$
by
\beq 
\vec w_{(i,j)} &=& \vec u_i \vec v_j,
\eeq
where the notation $(\vec v \tp \vec u)_{(i,j)}$ is used to
distinguish the 2-tensor $\vec w$ from an $M\times N$ matrix. The
tensor $\vec w$ is a vector in $\R^{MN}$ indexed by
pairs $(i,j)\in\nn{M}\times\nn{N}$; its norm is naturally defined to
be
\beqn
\label{eq:tensnorm}
\nrm{\vec w} = {\ts\oo2}\sum_{(i,j)\in\nn{M}\times\nn{N}}\abs{\vec w_{(i,j)}}.
\eeqn

To develop a convenient tensor notation, we will fix the index set
$V=\{1,\ldots,n\}$. 
For $I\subset V$, a tensor indexed by $I$ is a 
%linear map 
%$\vec u:\X^I\to\R$. 
vector $\vec u\in\R^{\X^I}$.
A special case of such an 
$I$-tensor is the product
$\vec u = \TP_{i\in I} \vec v\supr i$, where $\vec v\supr i\in\R^\X$ and
%
%More generally, if $I_1,I_2,\ldots I_K$ are 
%%finite 
%index sets and $\vec v\supr k\in\R^{I_K}$ for $1\leq k\leq K$, then
%the $K$-tensor $\vec w = \TP_{k=1}^K \vec v\supr k$ is defined by
\beq
%\vec w_{(i_1,i_2,\ldots,i_K)} &=& \prod_{k=1}^K \vec v\supr{k}_{i_k}
%[\vec u_E]_{(i_1,i_2,\ldots,i_{|E|})} &=&
%\prod_{j\in E} [\vec v\supr{j}]_{i_j}.
\vec u\bsub{x\rsub{I}} &=&
\prod_{i\in I} 
%(\vec v\supr{i})_{x_i}
\vec v\supr{i}\bsub{x_i}
\eeq
for each $x\rsub{I}\in{\X^I}$.
%with norm
%\beq
%\nrm{\vec w} = {\ts\oo2}\sum_{
%\abs{\vec w_{(i,j)}}.
%\eeq
To gain more familiarity with the notation, let us write the total
variation norm of an $I$-tensor:
\beqn
\label{eq:Itensnorm}
\nrm{\vec u} &=&
{\ts\oo2}\sum_{x\rsub{I}\in\X^I} \abs{\vec u\bsub{x\rsub{I}}}.
\eeqn
In order to extend
Lemma \ref{lem:tvtens}
to product tensors, we will need to define 
%a certain 
the
function
$\inex_k:\R^k\to\R$ and state some of its properties:
\belen
\label{lem:inex}
Define $\inex_k:\R^k\to\R$ recursively as $\inex_1(x)=x$ and
\beqn
\label{eq:inex}
\inex_{k+1}(x_1,x_2,\ldots,x_{k+1})=
x_{k+1}+(1-x_{k+1})\inex_{k}(x_1,x_2,\ldots,x_{k}).
\eeqn
Then
\bit
\item[\iia] $\inex_k$ is symmetric in its $k$
arguments, so it is well-defined as a mapping 
$$\inex:\{x_i:1\leq i\leq k\}\mapsto\R$$ from finite real sets to the reals
\item[\iib] $\inex_k$ takes $[0,1]^k$ to $[0,1]$ and is monotonically
      increasing in each argument on $[0,1]^k$
\item[\iic] If $B\subset C\subset[0,1]$ 
are finite sets
then $\inex(B)\leq\inex(C)$
\item[\iid] $\inex_k(x,x,\ldots,x)=1-(1-x)^k$
\item[\iie] if $B$ is finite and $1\in B\subset[0,1]$ then $\inex(B)=1$.
\item[\iif] if $B\subset[0,1]$
is a finite set then
 $\inex(B)\leq\sum_{x\in B} x$.
\eit
\enlen
\begin{rem}
In light of~\iia, we will use the notation
$\inex_k(x_1,x_2,\ldots,x_{k})$ and $\inex(\{x_i:1\leq i\leq k\})$
interchangeably, 
as dictated by convenience.
\end{rem}
\bepf
Claims~\iia,~\iib,~\iie,~\iif~are straightforward to verify from the
recursive 
definition of $\inex$ and induction. Claim~\iic~follows from~\iib~since
$$
\inex_{k+1}(x_1,x_2,\ldots,x_{k},0)=\inex_{k}(x_1,x_2,\ldots,x_{k})$$
and~\iid~is easily derived from the binomial expansion of $(1-x)^k$.
\enpf

The function $\inex_k$ is the natural generalization of 
$\inex_2(x_1,x_2)=x_1+x_2-x_1x_2$
to $k$ variables, and it is what we need 
for the analog of Lemma \ref{lem:tvtens}
for a  product of $k$ tensors:
\becon
\label{cor:tp}
%Suppose we have two sets of tensors 
Let
$\{\vec u\supr i\}
_{i\in I}
$ 
and 
$\{\vec v\supr i\}
_{i\in I}
$
be
two sets of tensors
%over the same index set
%$I=\{i_1,i_2,\ldots,i_k\}$
and
assume 
%futher 
that each of 
$\vec u\supr i,\vec v\supr i$ is a probability measure on $\X$. Then
we have
\beqn
\nrm{\TP_{i\in I}\vec u\supr i - \TP_{i\in I}\vec v\supr i}
&\leq&
\inex
%_{k}(
\tlprn{
\dsnrm{\vec u\supr i-\vec v\supr i}: i\in I
}
\hide{
\dsnrm{\vec u\supr1-\vec u\supr1},
\dsnrm{\vec u\supr2-\vec u\supr2},
\ldots,
\dsnrm{\vec u\supr{k}-\vec u\supr{k}})
}.
\eeqn
%where $\inex$ is being used in the sense of Lemma~\ref{lem:inex}\iia.
\hide{
where $\inex_k:\R^k\to\R$ is defined recursively as $\inex_1(x)=x$ and
\beqn
\label{eq:inex}
\inex_{k+1}(x_1,x_2,\ldots,x_{k+1})=
x_{k+1}+(1-x_{k+1})\inex_{k}(x_1,x_2,\ldots,x_{k}).
\eeqn
\encon
\begin{rem}
\label{rem:inex}
It is clear that the function $\inex_k$ is symmetric in its $k$
arguments, so it is well-defined as a mapping 
$$\inex:\{x_i:1\leq i\leq k\}\mapsto\R$$ from finite real sets to the reals.
\end{rem}
}
\encon
\bepf
Pick an $i_0\in I$ and let 
$\vec p=\vec u\supr{i_0}$,
$\vec q=\vec v\supr{i_0}$,
$$\vec p' = \TP_{i_0\neq i\in I} \vec u\supr i,
\qquad
\vec q' = \TP_{i_0\neq i\in I} \vec v\supr i.$$
Apply Lemma \ref{lem:tvtens} 
to 
$\nrm{\vec p\tp\vec q-\vec p'\tp\vec q'}$ and proceed by induction.
\enpf

Our final generalization concerns linear operators over
$I$-tensors. 
For $I,J\subseteq V$,
%An 
an
$I,J$-matrix $\vec A$ 
has dimensions $|\X^J|\times|\X^I|$ and
takes an $I$-tensor $\vec u$ to
a $J$-tensor $\vec v$:
for each $y\rsub{J}\in\X^J$,
we have
\beqn
\label{eq:Au}
\vec{v}\bsub{y\rsub{J}} 
&=&
\sum_{x\rsub{I}\in\X^I} 
\vec A\bsub{y\rsub{J},x\rsub{I}} \vec u\bsub{x\rsub{I}},
\eeqn
which we write as $\vec{Au}=\vec{v}$.
If $\vec A$ is an $I,J$-matrix and $\vec B$ is a $J,K$-matrix, the
matrix product $\vec{BA}$ is defined analogously to (\ref{eq:Au}).

As a special case, an $I,J$-matrix might factorize as a tensor product
of $|\X|\times|\X|$ matrices 
%$\vec A\supr{i,j}:\R^\X\to\R^\X$.
$\vec A\supr{i,j}\in\R^{\X\times\X}$.
We will write such a factorization in terms of a bipartite
graph\footnote{
Our notation for bipartite graphs is standard; it is equivalent to
$G=(I\cup J,E)$ where $I$ and $J$ are always assumed to be disjoint.
}
$G=(I+J,E)$, where $E\subset I\times J$ and the factors 
$\vec A\supr{i,j}$ are indexed by $(i,j)\in E$:
\beqn
\label{eq:matensprod}
 \vec A = \TP_{(i,j)\in E} \vec A\supr{i,j}, 
\eeqn
where
$$
%\vec A_{x_J,x_I} 
\vec A\bsub{y\rsub{J},x\rsub{I}}
= 
\prod_{(i,j)\in E} \vec A\supr{i,j}_{y_j,x_i}
$$
for all $x\rsub{I}\in\X^I$ and $y\rsub{J}\in\X^J$.
The norm of an $I,J$-matrix is a natural generalization of the matrix
norm defined in (\ref{eq:matnorm}):
\beqn
\label{eq:ijmatnorm}
\nrm{\vec A} = \max_{x\rsub{I},x\rsub{I}'\in\X^I} 
\nrm{
%\vec A_{*,x\rsub{I}}-\vec A_{*,x\rsub{I}'}
\vec A\bsub{*,x\rsub{I}} - 
\vec A\bsub{*,x\rsub{I}'}
}
\eeqn
where 
%$\vec A_{*,x_I}$ 
$\vec u=\vec A\bsub{*,x\rsub{I}}$
is the $J$-tensor given by
%$$ \vec u_{x_J} = \vec A_{x_J,x_I};$$
$$ \vec u\bsub{y\rsub{J}} = 
\vec A\bsub{y\rsub{J},x\rsub{I}};
$$
(\ref{eq:ijmatnorm}) is well-defined via the tensor norm in 
(\ref{eq:Itensnorm}).
Since $I,J$ matrices act on $I$-tensors by
ordinary matrix multiplication, 
%(\ref{eq:contr}) 
$\nrm{\vec{Au}} \leq \nrm{\vec A}\nrm{\vec u}$
continues to hold
when $\vec A$ is a
column-stochastic
$I,J$-matrix and $\vec u$ is a
balanced
$I$-tensor; if, additionally, $\vec B$ is a column-stochastic $J,K$-matrix,
%(\ref{eq:AB}) 
$\nrm{\vec{BA}} \leq \nrm{\vec B}\nrm{\vec A}$
also holds.
Likewise, since another way of writing
(\ref{eq:matensprod}) is
\beq
%\vec A_{*,x_I} = \TP_{(i,j)\in E} \vec A\supr{i,j}_{*,x_i},
\vec A\bsub{*,x\rsub{I}} = \TP_{(i,j)\in E}
\vec A\supr{i,j}
%\rsub{*,x_i}
\bsub{*,x_i}
,
\eeq
Corollary \ref{cor:tp} extends to tensor products of matrices:
\belen
\label{lem:TP}
Fix index sets 
$I,J$ 
%and 
%an edge set 
%$E=\{e_1,e_2,\ldots,e_k\}\subset I\times J$;
%these determine the
and a 
bipartite graph $(I+J,E)$. 
Let 
$\tlprn{\vec A\supr{i,j}}_{(i,j)\in E}$ be a collection of column-stochastic
$|\X|\times|\X|$ matrices, whose tensor product is the $I,J$ matrix
$$
 \vec A = 
\TP_{(i,j)\in E} \vec A\supr{i,j}.
%\TP_{t=1}^k \vec A\supr{i_t,j_t}.
$$
Then
\beq
\nrm{\vec A} &\leq&
\inex
\tlprn{
\dsnrm{\vec A\supr{i,j}}:(i,j)\in E
}
\hide{
(
\dsnrm{\vec A\supr{i_1,j_1}},
\dsnrm{\vec A\supr{i_2,j_2}},
\ldots,
\dsnrm{\vec A\supr{i_k,j_k}})
}.
\eeq
\enlen

We are now in a position to state the main technical lemma, from which
Theorem \ref{thm:mtmain} will follow straightforwardly:
\belen
\label{lem:maintech}
Let $\X$ be a finite set and
let $(X_i)_{1\leq i\leq n}$, $X_i\in\X$  be a Markov tree process, defined by a
tree $T=(V,E)$ and transition kernels $p_0$, 
$\tlprn{p_{uv}(\cdot\gn\cdot)}_{(u,v)\in E}$.
Let the $(u,v)$-contraction coefficient
$\tha_{uv}$
be as defined in (\ref{eq:uvthadef}).

Fix $1\leq i<j\leq n$ and let 
%$T_i=(V,E)$ be the subtree induced by
%$i$. 
%Then $T_i$ has root $i$ and 
%%suppose it is partitioned into $K$
%%levels: $V=\bigcup_{\el=0}^{K} \lev_{T_i}(\el)$. 
%levels $\{\lev_{T_i}(\el)\}_{\el=0,1,\ldots}$.
%Let 
$j_0=j_0(i,j)$ be as defined in Lemma \ref{lem:j0} (we are
assuming its existence, for otherwise $\etab_{ij}=0$).
Then we have
\beqn
\label{eq:etalpha}
\etab_{ij} &\leq&
\prod_{d=
%1
\depth_T(i)+1
}^{
\depth_{T}(j_0)
%d_0
}
\inex
\tlprn{
\tha_{uv} : v\in\lev_{T}(d)
}.
\eeqn
%where $\depth_{T}(\cdot)$ is defined in \S\ref{sec:graphtheo}.
%where $d_0=\depth_{T_i}(j_0)$.
\enlen
\newcommand{\sz}{\sseq{z}{i+1}{j-1}}
\newcommand{\sx}{\sseq{x}{j}{n}}
\bepf
For $y\in\X^{i-1}$ and $w,w'\in\X$, we have
\beqn
\eta_{ij}(y,w,w') &=&
{\ts\oo2}
\sumseq{x}{j}{n}
\abs{
\pr{\sseq{X}{j}{n}=\sseq{x}{j}{n}\gn\sseq{X}{1}{i}={y}{w}}
-    
\pr{\sseq{X}{j}{n}=\sseq{x}{j}{n}\gn\sseq{X}{1}{i}={y}{w'}}
}\\
%%%%%%%%%%%%%%%%%%%%%%%%%%%%%%%%
&=&
{\ts\oo2}
\sumseq{x}{j}{n}
\labs
\sumseq{z}{i+1}{j-1}
\lparen{
\pr{\sseq{X}{i+1}{n}=\sz \sx\gn\sseq{X}{1}{i}={y}{w}}
}\nonumber\\
&&
\qquad\qquad\quad
-    
\pr{\sseq{X}{i+1}{n}=\sz\,\sx\gn\sseq{X}{1}{i}={y}{w'}}
\rparen\rabs.
\eeqn
Let $T_i$ be the subtree induced by $i$ and
\beqn
\label{eq:ZC}
Z=T_i\cap\{i+1,\ldots,j_0-1\}
\qquad\text{and}\qquad
C = \{v\in T_i:(u,v)\in E, u<j_0, v\geq j_0\}.
\eeqn
Then by Lemma \ref{lem:j0} and the Markov property, we get
% this is numerically ``confirmed'' in
% etaijXdtree22(P,m,n,i,j,y,w1,w2)
\beqn
\nonumber
\eta_{ij}(y,w,w') &=&\\
&&
\nonumber
%\label{eq:hZC}
\hspace{-2cm}
{\ts\oo2}
\sum_{x\bsub{C}}
\labs
\sum_{x\bsub{Z}}
\lparen{
\pr{X\bsub{C\cup Z}=x\bsub{C\cup Z}\gn X_i=w}
}
-    
\pr{X\bsub{C\cup Z}=x\bsub{C\cup Z}\gn X_i=w'}
\rparen\rabs
\\\label{eq:hZC}&&
\eeqn
(the sum indexed by $\{j_0,\ldots,n\}\setminus C$ marginalizes out).

Define $D=\{d_k:k=0,\ldots,|D|\}$ with 
$d_0=\depth_T(i)$,
$d_{|D|}=\depth_T(j_0)$
and $d_{k+1}=d_k+1$ for $0\leq k<|D|$.
For $d\in D$, 
%For $\el=\depth_T(i),\ldots,\depth_{T}(j_0)$,
let $I_d=T_i\cap \lev_{T}(d)$ 
and $G_d=(I_{d-1}+I_d,E_d)$ be the
bipartite graph consisting of the nodes in $I_{d-1}$ and $I_{d}$,
and the edges in $E$ joining them (note that 
$I_{d_0}=\{i\}$).

For $(u,v)\in E$, let $\vec A\supr{u,v}$ be the $|\X|\times|\X|$
matrix given by
\beq
\vec A\supr{u,v}_{x,x'} = p_{uv}(x\gn x')
\eeq
and note that $\nrm{\vec A\supr{u,v}}=\tha_{uv}$.
Then by the Markov property, for each 
$z\bsub{I_{d}}\in\X^{I_d}$ and
$x\bsub{I_{d-1}}\in\X^{I_{d-1}}$, 
$d\in D\setminus\{d_0\}$,
we have
\beq
\pr{X\rsub{I_{d}}=z\rsub{I_{d}}\gn X\rsub{I_{d-1}}=x\rsub{I_{d-1}}}
&=&
\vec A\supr{d}\bsub{z\rsub{I_{d}},x\rsub{I_{d-1}}},
\eeq
where
\beq
\vec A\supr{d} &=&
\TP_{(u,v)\in E_d} \vec A\supr{u,v}.
\eeq
Likewise, for 
%$d>0$, 
$d\in D\setminus\{d_0\}$,
\beqn
\pr{X\rsub{I_{d}}=x\rsub{I_{d}}\gn X_i=w}
&=&
\sum_{x\rsub{I_1}'}\sum_{x\rsub{I_2}''}\cdots\sum_{x\rsub{I_{d-1}}\supr{d-1}}
\nonumber\\&&
\pr{X\rsub{I_{1}}=x\rsub{I_{1}}'\gn X_i=w}
\pr{X\rsub{I_{2}}=x\rsub{I_{2}}''\gn X\rsub{I_{1}}=x\rsub{I_{1}}'}
\cdots
\nonumber\\
&&\pr{X\rsub{I_{d}}=x\rsub{I_{d}}\gn X\rsub{I_{d-1}}=x\rsub{I_{d-1}}\supr{d-1}}
\nonumber\\
&=& 
\label{eq:AA}
(\vec A\supr d\vec A\supr{d-1}\cdots\vec A\supr{d_1})\bsub{x_{I_d},w}.
\eeqn
%Let $U=T_i\cap\{j_0,j_0+1,\ldots,n\}$
%and $Z=T_i\cap\{i+1,\ldots,j_0-1\}$.
Define the 
(balanced)
$I_{d_1}$-tensor
\beqn
\label{eq:hvdef}
\vec h=
%\vec A\supr{d_1}_{*,w} - \vec A\supr{d_1}_{*,w'},
\vec A\supr{d_1}\bsub{*,w} - \vec A\supr{d_1}\bsub{*,w'},
\eeqn
the $I_{d_{|D|}}$-tensor
\beqn
\label{eq:fvdef}
 \vec f = \vec A\supr{d_{|D|}}\vec A\supr{d_{|D|-1}}\cdots\vec
A\supr{d_2}\vec h,
\eeqn
and $C_0,C_1,Z_0
%,Z_1 
\subset 
\{1,\ldots,n\}
%\nn{n}
$:
\beqn
C_0=C\cap I_{\depth_T(j_0)},
\qquad
C_1 = C\setminus C_0,
\qquad
Z_0 = I_{\depth_T(j_0)}\setminus C_0,
%\qquad
%Z_1 = Z\setminus Z_0,
\eeqn
where $C$ and $Z$ are defined in (\ref{eq:ZC}).
%be the $I_1$-tensor defined by 
For readability we will write $\P(x_U\gn\cdot)$ instead of
$\pr{X_U=x_U\gn\cdot}$ below; no ambiguity should arise.
Combining (\ref{eq:hZC}) and (\ref{eq:AA}), we have
\beqn
\eta_{ij}(y,w,w') &=&
{\ts\oo2}
%\sum_{x_{C_0}}
%\sum_{x_{C_1}}
\sum_{x\rsub{C}}
\dsabs{
%\sum_{x_{Z_0}}
%\sum_{x_{Z_1}}
\sum_{x\rsub{Z}}
\paren{
\P(x\bsub{C\cup Z}\gn X_i=w)
-    
\P(x\bsub{C\cup Z}\gn X_i=w')
}}\\
&=&
{\ts\oo2}
\sum_{x\rsub{C_0}}
\sum_{x\rsub{C_1}}
\labs
\sum_{x\rsub{Z_0}}
\P(x\bsub{C_1}\gn x\bsub{Z_0})
\vec f\bsub{C_0\cup Z_0}
\rabs\\
&=&\nrm{\vec{Bf}}
\eeqn
where $\vec B$ is the
$|\X^{C_0\cup C_1}|\times|\X^{C_0\cup Z_0}|$ column-stochastic matrix
given by
$$ \vec B\bsub{x\rsub{C_0}\cup x\rsub{C_1},x'\rsub{C_0}\cup x_{Z_0}}
= \pred{x\rsub{C_0}=x'\rsub{C_0}} \P(x\rsub{C_1}\gn x\rsub{Z_0})$$
with the convention that $\P(x\rsub{C_1}\gn x\rsub{Z_0})=1$ if 
either of
$Z_0$ or $C_1$
is empty. The claim now follows by reading off the results previously
obtained:
\beq
\begin{array}{rcllll}
\nrm{\vec{Bf}}
&\leq& \nrm{\vec{B}}\nrm{\vec{f}}
&&&\text{Eq. (\ref{eq:contr})}\\\\
&\leq&\nrm{\vec{f}}
&&&\text{Remark \ref{rem:stochnorm}}\\\\
&\leq& \nrm{\vec h}\prod_{k=2}^{|D|} \nrm{\vec A\supr{d_k}}
&&&\text{Eqs. (\ref{eq:AB},\ref{eq:fvdef})}\\\\
&\leq& \prod_{k=1}^{|D|}
\inex\{\dsnrm{\vec A\supr{u,v}}:(u,v)\in E_{d_k}\}
&&&\text{Lemma \ref{lem:TP}}.
\end{array}
\eeq
\enpf

\bepf[Proof of Theorem \ref{thm:mtmain}]
We will borrow 
%all of 
the definitions from the proof of 
Lemma \ref{lem:maintech}.
To upper-bound $\etab_{ij}$ we first bound
$\inex\{\dsnrm{\vec A\supr{u,v}}:(u,v)\in E_{d_k}\}$. Since 
$$ |E_{d_k}|\leq \width(T)\leq L$$
(because every node in $I_{d_k}$ has exactly one parent in
$I_{d_{k-1}}$) and
$$ \nrm{\vec A\supr{u,v}}=\tha_{uv} \leq \tha < 1,$$ we appeal to
Lemma \ref{lem:inex} to obtain
\beqn
\inex\{\dsnrm{\vec A\supr{u,v}}:(u,v)\in E_{d_k}\}
&\leq&
1-(1-\tha)^L.
\eeqn
Now we must lower-bound the quantity
$h=\depth_T(j_0)-\depth_T(i)$. Since every level can have up to $L$ nodes,
we have
$$ j_0-i \leq hL $$
and so $h\geq\floor{(j_0-i)/L}\geq\floor{(j-i)/L}$.
\enpf

The calculations 
in Lemma \ref{lem:maintech}
yield
considerably more information than the simple bound in (\ref{eq:mainbd})
%in particular, the estimate in
---
certainly, the estimate in (\ref{eq:etalpha}) is quite a bit sharper.
Furthermore,
%For example, 
suppose the tree $T$ has 
%$K$ 
levels
%$\{I_d:d=0,\ldots,K\}$, 
$\{I_d:d=0,1,\ldots\}$
with the property that the levels are growing
at most linearly:
$$ |I_d|\leq cd $$
for some $c>0$. Let 
$d_i=\depth_T(i)$,
$d_j=\depth_T(j_0)$,
and $h=d_j-d_i$.
%$h=\depth_T(j_0)-\depth_T(i)$.
Then
\beq
j-i \leq j_0-i &\leq&
c\sum_{d_i+1}^{d_j} k \\
&=& \frac{c}{2}(d_j(d_j+1)-d_i(d_i+1))\\
&<& \frac{c}{2}((d_j+1)^2-d_i^2)\\
&<& \frac{c}{2}(d_i+h+1)^2
\eeq
so 
\hide{
There is no loss of generality in taking $i=1$ (since by
Lemma \ref{lem:j0}) only the nodes in $T_i$ affect the value of
$\eta_{ij}$. Then
$$ j-i \leq c\sum_{k=1}^h k = \frac{h(h+1)}{2} < c(h+1)^2/2, $$
so $$ h > \sqrt{2(j-i)/c}-1,$$
}
$$ h > \sqrt{2(j-i)/c}-d_i-1,$$
which yields the bound, via Lemma \ref{lem:inex}\iif,
\beqn
\etab_{ij} &\leq& \prod_{k=1}^h 
\sum_{(u,v)\in E_k} \tha_{uv}.
%(1-(1-\tha_k)^{ck})
\eeqn
Let
$\tha_k =\max\{\tha_{uv}: (u,v)\in E_k\}$;
then if
%If the contraction coefficients satisfy
$ck\tha_k\leq\beta$ 
holds
for some 
%$\tha<1$,
$\beta\in\R$,
%where $\tha_k \geq\max\{\tha_{uv}: (u,v)\in E_k\}$.
%When $\tha_k$ is small 
%($ck\tha_k\leq\tha< 1$), 
this becomes
\beqn
\nonumber
\etab_{ij} &
%\lesssim
\leq
& \prod_{k=1}^h (ck\tha_k)\\
\nonumber
&<&  \prod_{k=1}^{\sqrt{2(j-i)/c}-d_i-1} (ck\tha_k)\\
&\leq& 
\label{eq:thasqrt}
\beta^{\sqrt{2(j-i)/c}-d_i-1}.
\eeqn
%which 
This
is a non-trivial bound for trees with linearly growing
levels: recall that to bound $\nrm{\Delta}_\infty$,
%(\ref{eq:kontram},\ref{eq:infnorm}), 
%(\ref{eq:infnorm}), 
we must bound the series
$$ \sum_{j=i+1}^\infty \etab_{ij}.$$
By the limit comparison test with the series $\sum_{j=1}^\infty
1/j^2
$, we have that
$$ \sum_{j=i+1}^\infty \beta^{\sqrt{2(j-i)/c}-d_i-1}$$
converges for $\beta<1$.
Similar techniques may be applied when the level growth is
bounded by other slowly increasing functions.
It is hoped 
that 
%these techniques 
this method
will be extended to obtain concentration bounds
for larger classes of directed acyclic graphical models.

\section*{Acknowledgments}
Many
thanks to 
Gideon Schechtman for hosting and guidance at the Weizmann Institute,
and to Katalin Marton,
Elchanan Mossel, and Cosma Shalizi for helpful comments on the manuscript.
Thanks to Roi Weiss and an anonymous referee for insisting on a correct proof of Lemma \ref{lem:auv}.

\bibliographystyle{plain}
\bibliography{mark-conc}

\end{document}